\documentstyle[11pt]{article}
\begin{document}

\vskip 0.3cm
\centerline{\bf The duals and dual hulls}
\centerline{\bf of classes of analytic functions}

\vskip 0.3cm
\centerline{I.R.Nezhmetdinov}
\vskip 0.3cm

\centerline{1. Introduction}
\vskip 0.2cm

Let $D(a,R)= \{z:|z-a|<R \},$ where $a \in \mbox{\bf C}, R>0.$
In what follows we set for brevity: $D_R =D(0,R), D=D_1.$ 
Denote by ${\cal A}(D_R)$ the class of functions of the form 
$f(z)=\sum_{k=0}^{\infty} a_k (f) z^k,$
regular in $D.$  It is well known (see [1, Chapter 4, p. 98]) that
${\cal A}(D_R)$ is a topological vector space, the convergence in which 
is equivalent to locally uniform convergence in $D_R.$  
We also set ${\cal A}_0 (D_R) = \{ f \in {\cal A} (D):
a_0 (f)=1 \}, {\cal A} = {\cal A}(D), {\cal A}_0 ={\cal A}_0 (D),
{\cal A}(\overline{D})= \{ f \in {\cal A}: f \mbox{ is regular in 
the closed unit disk } \overline{D} \}.$

The space $\Lambda$ of the continuous linear functionals on 
${\cal A}$ is characterized by the following

Theorem A. (see [2]) A functional $\lambda$ on ${\cal A}$
is continuous and linear if and only if 
a function $g \in {\cal A}( \overline{D})$ can be found such that 
$\lambda(f)=(f*g)(1)$ for all $f \in {\cal A}.$ 
Here  
$(f*g)(z)= \sum_{k=0}^{\infty} a_k (f) a_k (g) z^k$ is the 
convolution (Hadamard product). In the sequel the correspondence 
between $\lambda$ and $g$ will be denoted as $\lambda := g.$

For $V \subset {\cal A}_0$ we define, in accordance with [3], 
the dual of $V$ as 
$V^* = \{ g \in {\cal A}_0: (f*g)(z) \neq 0, \quad
 \forall z \in D, \forall f \in V \}.$
We shall call $V$ a dual class if $V=W^*$ for some 
$W \subset {\cal A}_0.$ Various well known classes of functions 
can be represented in terms of the duality, thus new criteria for 
univalence, starlikeness, convexity etc.  arising. 
The dual hull of $V$ is the class $V^{**} = (V^* )^*$
which is the smallest of all dual classes containing $V.$ 
We state the duality principle (see [4]) which characterizes the dual hull 
for certain restrictions on $V.$

Theorem B. Let $V \subset {\cal A}_0$ be compact and satisfy the condition 
$$P_x f \in V \mbox{ for all } f \in V, \, x \in \overline{D},
\eqno{(1)}$$
where $(P_x f)(z) = f(xz), z \in D.$
Then for any $\lambda \in \Lambda$ we have 
$\lambda(V) = \lambda(V^{**}),$ moreover, 
$$f \in V \,  \Leftrightarrow \,  \forall \lambda \in  \Lambda
\quad \lambda(f) \in \lambda(V). \eqno{(2)} $$

In the present paper we obtain some new representations for duals 
and dual hulls. The duality principle is shown to be valid even under  
somewhat weakened conditions as compared with [3] and [4]. 
Finally, we introduce and consider subsets $U \subset V$ with 
$U^* =V^*$ (and $U^{**} = V^{**},$ resp.).

\vskip 0.3cm
\centerline{2. The principal results}
\vskip 0.2cm

In accordance with a definition from [4],  we shall say that 
aa set $V \subset {\cal A}_0,$  satisfying the condition (1), 
is complete. Define the complete hull of $V$ as a class 
$\mbox{cm}(V) = \cup_{x \in \overline{D}} P_x (V) =
\{ P_x f: f\in V, \, x\in \overline{D} \},$ being the smallest 
of all complete sets containing $V.$  It is easy to see that 
$ \mbox{cm} (V)$ coincides with the class $V',$
introduced in [3]. It was shown there that Theorem B holds true 
when (1) is replaced by a weaker condition 
$$\lambda (V)=\lambda [ \mbox{cm} (V)] \mbox{ for all }
\lambda \in \Lambda.  \eqno{(3)} $$

We state here without proof some elementary properties of the complete hull 
(for arbitrary sets $U, V \subset {\cal A}_0$):

a) if $V$ is compact, then $\mbox{cm}(V)$ is also compact;

b) $( \mbox{cm}(V))^* =V^*;$

c) $\mbox{cm}(U \cup V)=  \mbox{cm}(U) \cup \mbox{cm}(V);$

d) $\mbox{cm}(U \cap V) \subset \mbox{cm}(U) \cap \mbox{cm}(V).$

For $V \subset {\cal A}_0$ consider the class 
$V^T = \{ g \in {\cal A}_0 (\overline{D}): (f*g)(1) \neq 0
\mbox{ for all } f \in V \}.$  If $V$ is complete, then the completeness 
of $V^T$ easily follows but the converse is not true. In a similar way, 
for $U \subset {\cal A}_0 (\overline{D})$ introduce the class 
$U^{\perp} = \{ h \in {\cal A}_0 :(g*h)(1) \neq 0 
\mbox{ for all } g \in U \}.$

Now we state the principal results of the paper. 

\vskip0.1cm
Theorem 1. For $V \subset {\cal A}_0$ we have 
$V^* = \overline{(\mbox{cm} (V))^T },$
the closure being taken in the space ${\cal A}.$

\vskip0.1cm
Theorem 2. Let $V$ be a compact subclass of ${\cal A}_0,$ 
besides $V^T$ is complete. Then for any $\lambda \in \Lambda$ 
we have $\lambda(V) = \lambda(V^{**}),$ 
moreover, the equivalence (2) if valid.

\vskip0.1cm
Theorem 3. Under the conditions of the previous theorem there holds: 
$$ V^{**} = (V^T)^{\perp}. $$

\vskip 0.3cm
\centerline{3. The proof of the principal results}
\vskip 0.2cm

To begin with, let us prove several auxiliary assertions.

\vskip0.1cm
Lemma 1. If a sequence $f_n,$ $n=1,2,\dots,$ 
converges to $f$ in the space ${\cal A}(D_{R_1}),$ and  
$g_n \to g, \, n \to\infty$ in ${\cal A}(D_{R_2}),$ then  
$f_n * g_n \to f*g, $  $n \to\infty$ in ${\cal A}(D_{R_1 R_2} ).$

\vskip0.1cm
Proof. Assume that  
$z \in\overline{D_{\rho}},$ where  $0<\rho<R_1 R_2.$ Choose 
$\rho_1 < R_1, \rho_2 < R_2$ so that  
$\rho < \rho_1 \rho_2 <R_1 R_2.$ The convergence $f_n \to f$ 
in ${\cal A}(D_{R_1})$ implies the uniform convergence of 
$f_n$ to $f$ in $\overline{D_{\rho_1}},$ 
therefore, in its turn, the uniform boundedness of $\{ f_n \}$ follows in  
$\overline{D_{\rho_1}}.$ 
By virtue of the Cauchy inequalities for coefficients we can write 
$$ |a_k(f_n)| \le M(f_n ,\rho_1) \rho_1^{-k} \le M_1 \rho_1^{-k},     
\eqno{(4)} $$
$$|a_k(f_n) - a_k (f)| \le M(f_n -f ,\rho_1)\rho_1^{-k}  
\le \varepsilon_{1,n} \rho_1^{-k}  \mbox{ for all } n \ge 1, 
k \ge 0, \eqno{(5)}$$

\noindent 
where $M(f,r)= \sup_{|z|=r} |f(z)|, \varepsilon_{1,n} \to 0, n \to\infty.$ 
In the same way for the sequence $\{ g_n \}_{n=1}^{\infty}$ 
we get 

$$ |a_k (g)| \le M(g,\rho_2) \rho_2^{-k} \le M_2 \rho_2^{-k}, 
\eqno{(6)}$$
$$|a_k (g_n) - a_k (g)| \le M(g_n -g ,\rho_2)\rho_2^{-k} 
\le \varepsilon_{2,n} \rho_2^{-k} \mbox{ for all } n \ge 1,
k \ge 0,  \eqno{(7)}$$

\noindent 
where $\varepsilon_{2,n} \to 0, n \to\infty.$

By applying the estimates (4)-(7),  for all 
$z \in \overline{D_{\rho}}$  we have 

$$|(f_n *g_n )(z) - (f*g)(z)| \le 
\sum_{k=0}^{\infty} |a_k (f_n )a_k (g_n )-a_k(f)a_k(g)| \rho^k \le$$ 
$$\sum_{k=0}^{\infty} [|a_k (f_n)| |a_k (g_n )-a_k (g)| + 
|a_k (g)| |a_k (f_n) - a_k (f)|] \rho^k \le$$
$$\sum_{k=0}^{\infty} [\varepsilon_{2,n} \rho_2^{-k} M_1 \rho_1^{-k} 
+ \varepsilon_{1,n} \rho_1^{-k} M_2 \rho_2^{-k} ] \rho^k = 
(\varepsilon_{2,n} M_1  + \varepsilon_{1,n} M_2)(1- \rho/\rho_1 \rho_2 ).$$

The last expression tends to 0 as $n \to\infty,$ 
thus the proof is complete.

\vskip0.1cm
Lemma 2. Let $V$ be a compact subset of ${\cal A}_0(D_{R_1}), 
g \in {\cal A}(D_{R_2}).$
Then $U= \{ f*g: f \in V \}$ is compact in 
${\cal A}(D_{R_1 R_2}).$

\vskip0.1cm
Proof. Consider an arbitrary sequence of functions of the form 
$f_n *g, n=1,2,\dots,$ where $f_n \in V.$  By the compactness of   
$V$ a subsequence $f_{n_k}$ can be chosen to converge 
locally uniformly in $D_{R_1}$ to some function $f \in V.$ 
Then, by Lemma 1, $f_{n_k} * g \to f*g \in U$ 
in the space ${\cal A}(D_{R_1 R_2}),$ and, therefore, $U$ is compact.

\vskip0.1cm
Lemma 3. Let $V$ be a compact set in ${\cal A}(D_R), R>1,$ 
and in addition $f(1) \neq 0$ for any $f\in V.$  Then  
$\sigma \in (1,R)$ can be found such that $f(\sigma) \neq 0 $ 
for all $f \in V.$

\vskip0.1cm
Proof. Assuming the contrary, fix a decreasing number sequence 
$\{ x_n \}_{n=1}^{\infty}, 1 < x_n < R, n=1,2,\dots,$ converging to 1. 
Then for any $n \ge 1$ a function $f_n \in V$ can be found with 
$f_n (x_n )=0.$ Since $V$ is compact, then, choosing, if necessary, 
a subsequence, we may assume that $f_n \to f \in V$ in ${\cal A}(D_R).$  
Consider a sequence of functions $g_n (z)= (1-x_n z)^{-1},$ 
$n \in \mbox{\bf N},$ belonging to ${\cal A}(D_{1/x_1}).$  
It can be easily seen that $g_n (z) \rightarrow g(z)=(1-z)^{-1},$ 
$n \to \infty,$ locally uniformly in $D_{1/x_1},$ hence, 
by virtue of Lemma 1, $f_n *g_n \rightarrow f*g$ locally uniformly 
in the same disk. 
In particular, we have 
$$f_n (x_n )= (f_n *g_n )(1) \to (f*g)(1)= f(1),$$
whence $f(1)=0,$ thus contradicting to the conditions of the Lemma.

\vskip0.1cm
Lemma 4. Let $V \subset U \subset{\cal A}_0,$ besides for any 
$\lambda \in \Lambda, \lambda :=g,$ we have 
$g(0) \in \lambda (V).$ Then the following assertions are equivalent:

a) $\lambda (U)= \lambda (V)$ for all $\lambda \in \Lambda,$

b) if $\lambda \in  \Lambda, $then $0 \not\in \lambda (V)$ 
implies that $0 \not\in \lambda (U),$

c) $U^T =V^T.$

\vskip0.1cm
Proof. The implication a) $\Rightarrow$ b) is trivial. 
Assume that the condition b) holds. Then for $g \in V^T, 
\lambda:=g,$ we have $0 \not\in \lambda(V),$ 
therefore, $0\not\in\lambda(U),$ i.e.,  
$\lambda(f) = (g*f)(1) \neq 0$ for all $f \in U,$
thus $g\in U^T.$ Together with the obvious converse inclusion  
$U^T \subset V^T$ this yields $U^T = V^T.$ 

Finally, let c) hold. Fix a functional $\lambda \in  \Lambda, 
\lambda:=g.$ It suffices to show that  
$\lambda(U) \subset \lambda(V).$ Suppose that 
$w\in \mbox{\bf C} \setminus \lambda(V)$ (clearly, then  
$w \neq g(0)$). Hence for $f\in V$ we get  
$\lambda(f)-w = \{ [g(z)-w]*f(z) \}(1) \neq 0,$ 
consequently,
$$[g(z)-w][g(0)-w]^{-1} \in V^T = U^T. \eqno{(8)}$$

This is equivalent to the fact that $w \not\in \lambda(U),$ in other words,  
$\lambda(U) \subset \lambda(V),$ q.e.d.

Observe that the condition a) for $U = \mbox{cm}(V)$ coincides with 
(3). Show that for compact $V$ the additional assumption  
$g(0) \in  \lambda(V)$ may be discarded.

\vskip0.1cm
Lemma 5. Let $V \subset {\cal  A}_0$ be a compact class with 
$(\mbox{cm}(V))^T = V^T.$
Then $g(0) \in \lambda(V)$ for each $\lambda \in \Lambda, \lambda:=g.$

\vskip0.1cm
Let us prove the lemma by contradiction. Suppose that for some 
$\lambda \in \Lambda, \lambda:=g,$ we have $g(0) \not\in \lambda(V).$
By the compactness of $V$ and continuity of $\lambda$ the set $\lambda(V)$
is compact and a $\varepsilon>0$ can be found such that 
$D(g(0),\varepsilon) \cap \lambda(V) = \mbox{ \O}.$ 
Now, if $w \in D(g(0),\varepsilon), w \neq g(0),$ 
then, by reasoning as in the proof of Lemma 4, we deduce (8) 
with $U = \mbox{cm}(V).$ Thus, for all $f\in V, x \in \overline{D}$ 
we have 
$$\{ [g(z)-w][g(0)-w]^{-1} *P_x f \}(1)= 
[(g*f)(x)-w][g(0)-w]^{-1} \neq 0, $$
whence $(g*f)(x) \neq w$ for $x \in \overline{D}.$
Consequently, $g(0)=(g*f)(0)$ is an isolated point of the image $(g*f)(D)$ 
for any fixed function $f\in V.$ 
By the domain preservation principle, $(g*f)(z)$  is constant in $D$ 
and even in a somewhat larger disk. But then 
$\lambda(f) = (g*f)(1) = g(0),$
what contradicts to the assumption.

We remark that $V^T$ is complete if and only if 
$V^T = (\mbox{cm}(V))^T.$  Indeed, if $V^T$  is complete, then   
$P_x g \in V^T,$ provided that $g\in V^T, x \in \overline{D},$ 
hence, for every $f\in V$ we have 
$$(P_x g*f)(1)= (g*P_x f)(1) \neq 0, \eqno{(9)}$$
thus, $g\in (\mbox{cm}(V))^T.$ Since $V \subset \mbox{cm}(V),$ 
then by the above proved it follows that $V^T =(\mbox{cm}(V))^T.$ 
Conversely, if the latter equality holds, then for any 
$g\in V^T,$ in view of (9) with arbitrary $f\in V$  and 
$x \in \overline{D}$ we conclude that 
$P_x g \in V^T.$

\medskip
Proof of Theorem 1. Let $g\in V^*.$ Consider 
an increasing sequence $\{ r_n \}_{n=1}^{\infty}$ 
such that $ 0<r_n <1, n \in \mbox{\bf N},$  and  
$r_n \to 1, n \to \infty.$  
Put $g_n (z)= (P_{r_n} g)(z) = (1-r_n z)^{-1} * g(z).$ 
Since $0<r_n <1,$ we have  $g \in {\cal A}_0 (\overline{D}), 
n=1,2,\dots.$ 
It is easy to show that $(1-r_n z)^{-1} \to (1-z)^{-1}$ 
in ${\cal A},$ hence, by Lemma 1, we deduce that $g_n \to g, 
n \to\infty.$ If $f\in V,  x \in \overline{D},$  then
$(g_n *P_x f)(1)= (g*f)(r_n x)  \neq 0,$ because $g\in V^*.$
Therefore, $g_n \in (\mbox{cm}(V))^T,$  and 
$g\in  \overline{(\mbox{cm}(V))^T}.$ 

Conversely, let $g\in (\mbox{cm}(V)).$ Then for any 
$f\in V, x \in \overline{D}$ we have 
$(g*f)(x)= (g*P_x  f)(1)  \neq 0,$ i.e., $g\in V^*.$ 
The inclusion $(\mbox{cm}(V))^T  \subset V^*$  
in view of the closedness of $V^*$ in ${\cal A}$  (see [3]) 
implies that $\overline{(\mbox{cm}(V))^T}$  also lies in 
$V^*.$

Similar reasoning yield the following result, too.

\medskip
Theorem $1'.$ For $U\subset {\cal A}_0 \overline{(D)}$ 
we have $U^* = \overline{(\mbox{cm}(U))^{\perp}}.$

Observe that $U^T \subset U^{\perp},$ where $U^{\perp}$  
is not necessarily closed in ${\cal A}.$

\medskip
Proof of Theorem 2. By Lemmas 4 and 5 to prove the first assertion 
of the Theorem it suffices to verify that 
$V^T \subset (V^{**})^T.$ Let $g \in V^T = (\mbox{cm}(V))^T. $ 
Then $g\in {\cal A}(D_R), R>1,$ and for any function 
$f \in \mbox{cm}(V)$ we have $(g * f)(1) \neq   0.$ 
In view of the compactness of $\mbox{cm}(V)$ in ${\cal A},$ 
we conclude from Lemmas 2 and 3 that for certain 
$\sigma, 1< \sigma< R,$  the inequality $(g* f)(\sigma) \neq 0$ 
holds for all $f \in \mbox{cm}(V).$ 
Now, by setting $h=P_{\sigma} g,$ 
we have $h\in {\cal A} (\overline{D}),$ and 
$(h* f)(1) \neq 0,$ provided that $f \in  \mbox{cm}(V).$
Thus, $h \in (\mbox{cm}(V))^T,$ 
and, by Theorem 1, $h\in V^*.$ 
Hence, for an arbitrary $k \in V^{**}$ we have 
$$(g*k)(1) = (P_{1/\sigma} h*k)(1) = (h*k)(1/\sigma) \neq 0,$$
therefore $g \in (V^{**})^T.$

We prove now the second part of the Theorem. As it was proved earlier, 
for all $\lambda \in \Lambda$ we have $\lambda(V)= \lambda(V),$ 
and $f\in V^{**}$ implies $\lambda(f)\in  \lambda(V).$
Conversely, let for a function $f\in {\cal A}_0$ 
there holds $\lambda(f) \in \lambda(V)$ for all $\lambda \in  
\Lambda.$ Fix an arbitrary $g \in V^T, \lambda := g.$  
By the completeness of $V^T$ we have 
$\lambda(P_x h)= (g*P_x h)(1) \neq  0,$ 
whatever $h \in V, x \in \overline{D}.$
Put $\lambda_x  :=P_x g.$ For $h \in V, x \in \overline{D}$ 
we get 
$\lambda_x (h)=(P_x g*h)(1)= (g*P_x h)(1)  \neq  0,$ i.e., 
$0 \not\in \lambda_x (V).$ 
But then 
$\lambda_x (f) = (g*f)(x) \neq  0,$ and $f \in (V^T)^*  = 
[(\mbox{cm}(V))^T]^*  = V^{**}.$ 

The last equality may be proved if we show that 
$(\overline{U})^* =U^*$  for $U \subset {\cal A}_0.$  
The inclusion $U \subset \overline{U}$ implies that  
$(\overline{U})^* \subset U^*.$ 
Let us prove the converse. For any $f\in \overline{U}$ there exists 
a sequence $f_n \in U, n \in \mbox{\bf N},$ converging to 
$f.$  Now, if $g\in U^*,$ then  $(g*f_n )(z) \neq  0$ for 
$n\in \mbox{\bf N}, z\in D.$ By Lemma 1, $g*f_n \to g*f, 
n \to \infty,$ in ${\cal A}.$ According to Hurwitz theorem 
(see [5, p.19]), if $(g*f)(z) \not\equiv \, \mbox{const},$ 
then $(g*f)(z) \neq 0$ in $D,$ however, in the case when 
$(g*f)(z) \equiv \mbox{const},$ we also obtain 
$(g*f)(z) = (g*f)(0)= 1 \neq 0.$ 
Since the function $f\in \overline{U}$ is chosen arbitrarily, 
then $g\in (\overline{U})^*,$ and the proof is complete.

\medskip
Proof of Theorem 3. Let 
$h\in V ^{**}, g \in V^T, \lambda:=g.$ 
Then $\lambda(f)= (g*f)(1) \neq 0 $ for all $f\in V,$ i.e., 
$0 \not\in \lambda(V).$ By Theorem 2, $0 \neq \lambda(h)= 
(g*h)(1),$ whence $h\in (V^T)^{\perp}.$ 
On the other hand, let $h\in (V^T)^{\perp},$ and $g\in V^*.$ 
By virtue of Theorem 1 a sequence of functions 
$g_n \in V^T =(\mbox{cm}(V))^T, n=1,2\dots,$ can be found such that 
$g_n \to g, n \to\infty.$ In addition for all $x \in \overline{D}$  
we have:  $(g_n *P_x h)(1)= (g_n *h)(x) \neq 0.$ 
By using the Hurwitz theorem as in the  previous proof, we obtain  
$(h*g)(x) \neq 0,$ for all $x\in D,$ thus, $h\in V^{**}.$ 

\vskip0.1cm
Corollary 1. Under the assumptions of Theorem 2 the following 
equalities hold:
$$V^{**} = \cap_{\lambda \in \Lambda} \lambda^{-1} [\lambda(V)]= 
\cap_{\lambda \in \Lambda} (V + \, \mbox{ker} \lambda),$$
where $\lambda^{-1}(A)$ is the inverse image of the set $A$ 
for the mapping $\lambda, A+B= \{ x+y: \, x\in A, y\in B \}$ 
is an algebraic sum of two sets, and 
$\mbox{ker}\lambda= \{ f\in {\cal A}: \lambda(f)=0 \}$ 
is the kernel of the functional $\lambda.$ 

\vskip0.1cm
Proof. According to Theorem 2, the inclusion $f\in V^{**}$ 
is equivalent to the fact that $\lambda(f)\in \lambda(V)$ for any 
$\lambda \in \Lambda$ which, in its turn, equivalent to the 
inclusion 
$f\in \lambda^{-1} [\lambda(V)], \lambda \in \Lambda,$ 
and the first equality follows. Now, if 
$f\in \lambda^{-1} [\lambda(V)],$ then $g\in V$ 
can be found such that $\lambda(g)= \lambda(f).$ 
But then for $h=f-g$ we have 
$\lambda(h)= \lambda(f)- \lambda(g) =0,$  and
$h\in \mbox{ker}\lambda,$ whence 
$f= g+h \in V + \mbox{ker}\lambda.$ 
On the other hand, if $f=g+h,$ where 
$g\in V, h\in \mbox{ker}\lambda,$ then  
$\lambda(f)=\lambda(g)\in \lambda(V).$

Remark. If $V \subset {\cal A}_0$ is a compact dual class, then 
$V = \cap_{\lambda \in \Lambda} \lambda^{-1} [\lambda(V)].$
A similar representation is valid for an arbitrary compact 
convex set in locally convex space $X$, however, $\Lambda$ 
should be replaced here by the space of all 
real-valued continuous linear functionals on 
$X$ (see [6 , Chapter 2, p. 88]).

\medskip
Corollary 2. If $U, V \subset {\cal A}_0$ are compact classes and 
$U^T, V^T$ are complete, then the following relations are 
equivalent:

a) $U^{**} = V^{**};$ b) $U^T = V^T;$  c) $U^{*} = V^{*}.$

\medskip
Proof. Suppose that a) holds. Let $g \in U^T , \lambda:= g.$ 
Then, in view of the duality principle we have 
$0 \not\in \lambda(U)= \lambda(U^{**})= \lambda(V^{**})= 
\lambda(V),$ therefore, $g\in V^T.$ 
The converse inclusion is proved in a similar way. 
By virtue of Theorem 1 and the remark about completeness of  
$U^T$  we deduce that b) $\Rightarrow$ c). 
To conclude with, clearly, c) implies a).

\medskip
Theorem 4. Let $V\subset{\cal A}_0$ be a compact set and 
$V^T$ be complete. Then 
$$V = \cup_{0<r<1} P_r (V^T)= 
\cup_{0<r<1} P_r \{ [P_r (V)]^T\}. \eqno{(10)}$$

\medskip
Proof. By the completeness of $V^T$ it follows that for any 
$r\in (0,1)$ we have $P_r (V^T )\subset V^T,$ thus,
$\cup_{0<r<1} P_r (V^T ) \subset V^T.$ 
Let $g\in V^T.$ Then $g\in  {\cal  A}(D_R  )$  for certain 
$R>1,$ besides $(g*f)(1) \neq 0$ for all $f\in V.$ 
By Lemmas 2 and 3 $\sigma \in (1,R)$ can be found such that 
$(g*f)(\sigma) \neq 0  \quad \forall f\in V.$ 
Setting $h= P_{\sigma} g$ (here $h\in {\cal A}(\overline{D})),$ 
we have 
$(h*f)(1)=(g*f)(\sigma) \neq 0,$ if  $f\in V,$ and, thus, 
$h\in V^T.$ Therefore, $g = P_{1/\sigma} h 
\in P_{1/\sigma} (V^T) \subset \cup_{0<r<1} P_r (V^T).$

Let us prove the second equality. Fix $g\in V^T$ and $r \in (0,1).$ 
Since $V^T$ is complete, then for all $f \in V$
we have $(g*P_r f)(1) \neq 0,$ and, hence, $g\in [P_r (V)]^T.$ 
Now the inclusion $V^T \subset [P_r (V)]^T$ obviously yields 
$P_r (V^T ) \subset P_r \{ [P_r (V)]^T \},$ 
whence 
$$V^T = \cup_{0<r<1} P_r (V^T ) \subset 
\cup_{0<r<1} P_r \{ [P_r (V)]^T \}.$$
On the other hand, if $g$ belongs to the right-hand side of (10), 
then for  some $r\in (0,1)$ we get $g \in P_r \{ [P_r (V)]^T \},$ 
hence, $g=P_r h,$ where $h \in  [P_r (V)]^T.$ 
Clearly, $g \in {\cal A}_0(\overline{D}),$ and 
$(h*P_r f)(1) = (P_r h*f)(1) \neq 0,$ 
consequently, $g\in V^T,$ q.e.d.

\vskip 0.3cm
\centerline{4. Border elements of the class $V$}
\vskip 0.2cm

As it was said above, $V^{**}$ is the least (by inclusion)
dual class containing $V.$ Many extremal problems 
on various classes of analytic functions are solved by reducing 
to more elementary extremal problems on simpler subclasses, 
having the same closed convex hull as the given class. 
We construct a subclass $U \subset V$ with $U^* =V^*$  
and study some of its properties.

Let $V\subset{\cal A}_0$ and is distinct from the class consisting of 
the unique element $e \equiv 1.$ We shall call $f\in V$ 
a border element of $V,$ if the relation $f=P_x g$ 
with $g\in V, x \in \overline{D},$ implies that $|x|=1.$ 
The set of all border elements of $V$ will be called the 
border of $V$ and denoted by $\mbox{bor}(V).$ 
If $V= \{ e \},$ then, by definition, we set $\mbox{bor}(V)= \{ e \}.$

\medskip
Lemma 6. Let $V \subset {\cal A}_0$  be a compact set.
Then for any function $f\in V$ we can choose $g\in \mbox{bor}(V)$ 
$and x \in \overline{D}$ such that $f=P_x g.$ 
In addition, for $f \not\equiv e$ the element $g$ is determined 
modulo transformation $P_y$ with $|y|=1.$

\medskip
Proof. Fix any function $f\in V$ and consider a number set 
$$R_f = \{ r>0: \mbox{ there exist } g\in V, x \in \overline{D} 
\mbox{ such that } f=P_x g, |x|=r \}.$$ 
Note, that $1 \in R_f \neq \mbox{\O}.$ 
Let $r_0 = \inf R_f > 0.$ 
Then a sequence $\{ r_n \}_{n=1}^{\infty}$ 
of elements from $R_f$ can be found that converges to $r_0$ as $n\to\infty,$ 
consequently, sequences $g_n \in V, x_n \in D$ can be indicated 
such that $f = P_{x_n} g_n, |x_n| = r_ , n=1,2,\dots.$ 
By the compactness of the sets $V$ and $\overline{D}$ assume that 
(choosing, if necessary, subsequences) 
$g_n \to g_0 \in V, x_n \to x_0  \in \overline{D}, n\to\infty.$ 
It is easy to see that $r_0 =|x_0|,$ and, by Lemma 1, 
we have $P_{x_n} g_n  \to P_{x_0} g_0, n\to\infty,$ therefore,  
$f=P_{x_0} g_0.$ 
Suppose that $g \not\in \mbox{bor}(V).$ 
Then for certain $h\in V$ and $y\in D$ we get: $g_0  = P_y h,$ 
hence $f = P_{x_0 y} h,$ where $|x_0 y| \in R_f,$ 
however, $|x_0 y|< |x_0 |= r_0,$ thus contradicting the assumption.
Consequently, $g_0$  is the required element. 

It remains to consider the case $r_0= 0.$ 
By reasoning as above, we get $f = P_0 g_0  = e.$ 
If $V \neq \{ e \},$ then a function $h \in V,$ $h \neq e,$ 
can be found represented in the form $h = P_x g$ 
for some $g \in \mbox{bor}(V)$ and $x \in \overline{D}.$ 
But then $f= P_0 g,$ q.e.d. For the class $V= \{ e \}$ 
the reasoning is  trivial.

Let us prove the uniqueness of the required function. 
Let for some $g, h\in  \mbox{bor}(V)$ and $x,y\in \overline{D}$ 
we have $f = P_x g = P_y h,$ besides $f \neq e.$  
Clearly, $x,y \neq  0.$ Assume, for the definiteness' sake, that 
$|x| \le |y|.$  Since the equalities 
$$f(z)= g(xz)= h(yz)= g[(x/y)yz]= (P_{x/y} g)(yz), z \in D,$$
hold, then in the disk $D_{|y|}$ the functions $h$ and $P_{x/y} g$ 
coincide. By the uniqueness theorem they coincide in the whole 
disk $D.$ 
But then, since $h \in \mbox{bor}(V),$ we have $|x/y| = 1.$ 
The case when $|x| \ge |y|$ is studied similarly.

Remark. Lemma 6 implies that for any compact 
$V\subset{\cal A}_0$ the imclusion $V \subset \mbox{cm}(\mbox{bor}(V))$ 
holds, the equality being attained whenever $V$ is complete. 

\medskip
Corollary 3. If $V\subset{\cal A}_0$ is compact, then 
$(\mbox{bor}(V))^* =V^*.$

\medskip
Proof. By the definition and the immediately preceding remark we have: 
$\mbox{bor}(V) \subset V \subset \mbox{cm}(\mbox{bor}(V)),$  
therefore, $[\mbox{cm}(\mbox{bor}(V))]^*  \subset V^* \subset 
(\mbox{bor}(V))^*.$ 
In view of the property b) of the complete hull, we deduce that 
the exterior terms of the last inclusion coincide, q.e.d. 

\medskip
Theorem 5. Let $V\subset{\cal A}_0$ be compact and let $V^T$ 
be a complete set. Then for any $\lambda \in \Lambda, 
\lambda:=g,$ the relations hold: 
$$\lambda(V)= \cup_{f \in \mbox{bor}(V)} (f*g)(\overline{D}), 
\eqno{(11)}$$
$$\partial\lambda(V) \subset \cup_{f \in \mbox{bor}(V)}
(f*g)(\partial D), \eqno{(12)} $$

\medskip
Proof. Suppose that $c \in \lambda(V),$ where 
$\lambda \in \Lambda, \lambda:=g.$  It means that $f\in V$  
can be found such that $c=\lambda(f)= (g*f)(1).$  
By Lemma 6, the function $f$ can be presented in the form $f=P_x h,$ 
where $h\in \mbox{bor}(V), x \in \overline{D}.$ 
But then we get 
$$c= (f*g)(1)= (P_x h*g)(1)= (h*g)(x)\in (h*g)(\overline{D}), 
\eqno{(13)}$$
where the last set, clearly, contains in the right-hand side of (11).

Conversely, let 
$c \in \cup_{f \in \mbox{bor}(V)} (f*g)(\overline{D}).$ 
Then there exist 
$f\in \mbox{bor}(V) \subset V$ and $x \in \overline{D}$ 
such that 
$$c= (f*g)(x)= (P_x f*g)(1)= \lambda(P_x f). $$

Since $P_x f \in \mbox{cm}(V) \subset V,$ 
then, by Theorem 2, $c= \lambda(P_x f)$ lies in 
$\lambda(V)$, q.e.d.

We prove now the inclusion (12). Assume that 
$c \in \partial\lambda(V),$ then $c\in \lambda(V),$ 
because $\lambda(V)$ is compact. In view of (11),
$f \in  \mbox{bor}(V), x \in \overline{D},$ can be found to 
fulfil the  relations (13).
Here, if $f*g \not\equiv \mbox{const}$ and 
$x \in D,$ then, by the domain conservation principle, 
certain neighborhood $O(c)$ of the point $c$ lies in $(f*g)(D),$  
and for any $c'\in O(c)$ we can find $x'\in D$ such that 
$c'= (f*g)(x')= \lambda(P_{x'} h).$ 
Hence, as above, $\lambda(P_{x'} f) \in \lambda(V),$ 
i.e., $c$ is an interior point of $\lambda(V),$ 
what contradicts to the assumption. Therefore, $x$ 
must be a boundary point of $D.$ If $f*g = \mbox{const},$ 
then $c = (f*g)(1) \in (f*g)(\partial D).$

An example of the compact class $V = \{ 1 +xz: x \in \overline{D} \} \cup 
\{ 1 +yz^2: y \in \overline{D} \}$ and the functional 
$\lambda (f) = a_1 (f), \lambda := g(z) = z,$ 
shows that the converse inclusion in (12) does not hold. 

\medskip
Theorem 6. For any compact $V\subset{\cal A}_0$ 
we have 
$$\mbox{bor}(V^* )= V^* \setminus ( \mbox{cm}(V))^T.$$

\medskip
Proof. From Theorem 1 and the definition of the border it follows that 
both $(\mbox{cm}(V))^T$ and $\mbox{bor}(V^* )$ 
are subsets of $V^*.$ 
Let $g \in (\mbox{cm}(V))^T.$ Applying Lemmas 2 and 3, as in the 
proof of Theorem 2, we conclude that for some $\sigma >1$ 
there holds $h= P_{\sigma} g \in V^*,$
hence $g = P_{1/\sigma} h$, and therefore 
$g \not\in \mbox{bor}(V^*).$ 

Conversely, if $g \in  V^* \setminus \mbox{bor}(V^* ),$ 
then choose $ h\in V^*$ and $x\in D$  such that 
$g=P_x h.$ Then $g\in {\cal A}_0 (\overline{D}),$ 
and simultaneously for all $y\in \overline{D}, f\in V$ 
we have 
$$(g*P_y f)(1) = (P_x h*P_y f)(1) = (h*f)(xy)  \neq  0, 
\eqno{(14)} $$
because $h\in V^*.$ But (14) implies that 
$g \in (\mbox{cm}(V)),$ q.e.d.

\vskip 0.3cm
\centerline{References}
\vskip 0.3cm

1. Hallenbeck D.J., MacGregor T.H. Linear problems and
convexity technique. - New York: Pitman Publishers, 1984, 182 pp.

2. Toeplitz O. Die linearen vollkommenen Ra\"ume der 
Funktionentheo\-rie //Comment. Math. Helv. - 1949. - V.23. - 
P.222--242.

3. Ruscheweyh St. Duality for Hadamard products with 
applications to extremal problems for functions regular  in  
the unit disc//Trans. Amer. Math. Soc. - 1975. - V.210 - P.63--74.

4. Ruscheweyh St. Convolutions in geometric function theory.
- Montr\'eal: Les Presses de l'Universit\'e de Montr\'eal, 1982.
- 166 pp.

5. Golusin G.M. Geometric theory of functions of complex variables. 
- 2nd ed. - Moscow.: Nauka, 1966. - 628 p.

6. Schaefer H. Topological linear spaces. - Moscow: Mir, 1971. 
- 359 p.

\end{document}